\documentclass[b5paper,twoside]{article}
\usepackage{amsmath, amsthm, amssymb,pstricks,pst-node,pst-plot}
\usepackage[dvips]{graphicx}
\def\titlePS#1{\noindent{{\LARGE\sffamily #1}} \vspace{3mm}}
\def\shorttitlePS#1{\markboth{Prague Stochastics 2006}{#1}}
\def\amsPS#1{\noindent{\it MSC 2000:} #1}
\def\keywPS#1{\noindent{\it Key words:} #1 \\}

\def\thanksPS#1{\footnotetext{#1}}
\def\authorsPS{}
\newcommand\authPS[8]{\ifnum\authPSc<1\def\authPSc{1}\else, \fi{\large\sffamily #1 #2}%
\edef\authorsPS{\authorsPS \par\vskip 1mm\noindent #1 #2:\hskip 5mm #3, #4, #5, #6, #7, #8}}
\newenvironment*{abstractPS}{\vspace{3mm}\noindent{\it Abstract:} }{\vspace{3mm}}
\newtheorem{theorem}{Theorem}[section]

\newtheorem{proposition}[theorem]{Proposition}
\theoremstyle{definition}
\newtheorem{definition}[theorem]{Definition}
\theoremstyle{remark}

\newtheorem{example}[theorem]{Example}
\def\authPSc{0}
\begin{document}


\titlePS{Algebraic Techniques for Gaussian Models}

\shorttitlePS{Algebraic Techniques for Gaussian Models}

\thanksPS{Research supported by the US National Science Foundation (DMS-0505612).}

\noindent
\authPS{Mathias}
{Drton}
{The University of Chicago, Department of Statistics}
{5734
  S.~University Ave}
{Chicago, Illinois}
{60637}
{U.S.A.}
{drton@galton.uchicago.edu}%


\begin{abstractPS}
  Many statistical models are algebraic in that they are defined by
  polynomial constraints or by parameterizations that are polynomial or
  rational maps.  This opens the door for tools from computational
  algebraic geometry.  These tools can be employed to solve equation
  systems arising in maximum likelihood estimation and parameter
  identification, but they also permit to study model singularities at
  which standard asymptotic approximations to the distribution of
  estimators and test statistics may no longer be valid.  This paper
  demonstrates such applications of algebraic geometry in selected examples
  of Gaussian models, thereby complementing the existing literature on
  models for discrete variables.
\end{abstractPS}

\amsPS{62H05, 62H12}

\keywPS{Algebraic statistics, multivariate normal distribution, parameter
  identification, singularities}


\newcommand{\ind}{\mbox{$\perp \kern-5.5pt \perp$}}
\newcommand{\nind}{\mbox{$\not\hspace{-4pt}\ind$}}

\section{Introduction}
\label{sec:intro}

Algebraic statistics applies algebraic geometry to gain insight in
structure and properties of statistical models, and to tackle computational
problems arising in tasks of statistical inference.  Work in this field has
addressed, for example, exact tests in contingency tables, experimental
design, phylogenetic trees, maximum likelihood estimation under multinomial
sampling, and Bayesian networks; cf.~\cite{ascb,pistone:01}.  Algebraic
geometry typically enters the playing field in one of two ways.  On one
hand, statistical models are sometimes derived from a simple saturated
model by imposing constraints.  These constraints may, in particular, be
motivated by considerations of (conditional) independence, stationarity or
homogeneity.  If the constraints are polynomial constraints on the
parameters of the saturated model, then the model corresponds to the
intersection of an {\em algebraic variety} and the saturated parameter
space.  An algebraic variety is the solution set of a system of polynomial
equations. On the other hand, many statistical models are defined via a
parameterization rather than via constraints.  However, if this
parameterization is a polynomial, or more generally a rational map, then
the model, which can be identified with the image of the parameterization
map, is naturally embedded in an algebraic variety.  This algebraic
description of the model is often useful because it can reveal insights
about the model that are not as readily obtained from the parameterization
alone.

Virtually all work in algebraic statistics considers purely discrete
variables; see \cite{algfac} for an exception.  The sample space is
then finite, and the objects of interest are algebraic varieties over
a probability simplex.
However, the philosophy of algebraic statistics applies, regardless of the
distributional setting, whenever a statistical model is an ``algebraic''
submodel of some natural supermodel.  A particularly interesting case
occurs when the supermodel is a regular exponential family, because 
algebraic submodels may inherit desirable statistical properties at points
at which the submodel's local geometry is sufficiently ``regular.''

Discussing simple problems from parameter identification and likelihood
ratio testing, this paper demonstrates algebraic techniques for Gaussian
models, i.e., families of (non-singular) multivariate normal distributions.
Section \ref{sec:multnorm} reviews the normal distribution and introduces
the algebraic point of view.  Section \ref{sec:ident} treats the problem of
identification of a graphical model with hidden variables.  Section
\ref{sec:sing} is devoted to model singularities, at which standard
$\chi^2$-approximations to the distribution of the likelihood ratio test
statistic may no longer be valid.

\section{Algebraic Gaussian models}
\label{sec:multnorm}

Let $\mathbb{R}^{p\times p}_{\rm pd}$ and $\mathbb{R}^{p\times p}_{\rm
  psd}$ be the cones of positive definite and positive semi-definite
symmetric $p\times p$-matrices, respectively.  The {\em multivariate normal
  distribution} $\mathcal{N}_p(\mu,\Sigma)$ with {\em mean vector}
$\mu=(\mu_1,\dots,\mu_p)^t\in\mathbb{R}^p$ and {\em covariance matrix}
$\Sigma=(\sigma_{ij})\in\mathbb{R}^{p\times p}_{\rm pd}$ is the probability
distribution on $\mathbb{R}^p$ that has Lebesgue density function
\begin{align*}
f_{\mu,\Sigma}(x) =
\frac{1}{\sqrt{(2\pi)^p\det(\Sigma)}} \;\exp\left\{
  -\frac{1}{2}(x-\mu)^t\Sigma^{-1}(x-\mu)\right\}, \qquad x\in\mathbb{R}^p.
\end{align*}
 A {\em Gaussian
  (statistical) model} with {\em mean parameter space\/} $M\subseteq
\mathbb{R}^p\times \mathbb{R}^{p\times p}_{\rm pd}$ is the family of
multivariate normal distributions $\left\{ \mathcal{N}_p(\mu,\Sigma)\mid
  (\mu,\Sigma)\in M \right\}$.

\begin{proposition}[\mbox{\cite[p.~194]{twa}}]
  \label{prop:expfam}
  The {saturated Gaussian model}, that is, the family of all
  multivariate normal distributions on $\mathbb{R}^p$, which has mean
  parameter space $M=\mathbb{R}^p\times \mathbb{R}^{p\times p}_{\rm
    pd}$, forms a regular exponential family with {sufficient
    statistics} $x\in\mathbb{R}^p$ and $xx^t\in \mathbb{R}^{p\times
    p}_{\rm psd}$. The {natural parameters} are
  $\Sigma^{-1}\mu\in\mathbb{R}^p$ and $\Sigma^{-1}\in
  \mathbb{R}^{p\times p}_{\rm pd}$.
\end{proposition}

Statistical modelling in the Gaussian framework involves hypotheses
about structural relationships among the components of the mean
parameters $\mu$ and $\Sigma$.  In many interesting cases, such a
relationship comes from a parameterization.  If this parameterization
is rational, as detailed in the following definition, then the
resulting model can be studied taking an algebraic point of view.
Recall that a set $\Theta\subseteq\mathbb{R}^d$ is semi-algebraic if
it is a union of sets of points satisfying polynomial equalities and
inequalities; compare Chapter 2 in \cite{benedetti}.

\begin{definition}
  A Gaussian model is a {\em parametric algebraic model} if its mean
  parameter space $M=\mathbf{f}(\Theta)$, where
  $\Theta\subseteq\mathbb{R}^d$ is an open semi-algebraic set, and
  \begin{align*}
    \mathbf{f}:\Theta &\to \mathbb{R}^p\times
    \mathbb{R}^{p\times p}_{\rm pd}\\
    \theta &\mapsto
    \left(\frac{g_1}{h_1}(\theta),\dots,\frac{g_p}{h_p}(\theta),\;
      \frac{g_{11}}{h_{11}}(\theta),\dots,\frac{g_{pp}}{h_{pp}}(\theta)\right)
  \end{align*}
  is a rational map defined everywhere on $\Theta$.  In other words, the
  functions $g_k$, $h_k$, $g_{ij}$ and $h_{ij}$ are polynomial functions
  such that $0\not\in h_k(\Theta)$ for all $k\in [p]:=\{1,\dots,p\}$ and
  $0\not\in h_{ij}(\Theta)$ for all $(i,j)\in [p]^2$.
\end{definition}

Not all statistical models of interest are specified in terms of a
parameterization; instead the model may be specified implicitly in the form
of constraints on the mean parameters $\mu$ and $\Sigma$.  One important
class of constraints arises from conditional independence, which in the
multivariate normal distribution corresponds to well-known polynomial
conditions on the covariance matrix $\Sigma$.

\begin{proposition}
  \label{prop:ci}
  Let $X$ be a random vector in $\mathbb{R}^p$ that follows a multivariate
  normal distribution $\mathcal{N}_p(\mu,\Sigma)$, in symbols,
  $X\sim\mathcal{N}_p(\mu,\Sigma)$.  For three pairwise disjoint index sets
  $A,B,C\subseteq [p]:=\{1,\dots,p\}$, it holds that
  \[
  X_A\ind X_B\mid X_C \quad \iff \quad \det(\Sigma_{\{i\}\cup C\times
    \{j\} \cup C}) = 0 \quad\forall i\in A,\; j\in B.
  \]
  Here, $X_A\ind X_B\mid X_C$ denotes conditional independence of $X_A$ and
  $X_B$ given $X_C$, and $X_A\ind X_B\mid X_\emptyset$ denotes marginal
  independence of $X_A$ and $X_B$.   
\end{proposition}

The fact that conditional independence is an algebraic condition
motivates the next definition, in which
\[
\mathbb{R}[\mu_k,\sigma_{ij}\mid i,j,k\in[p],\; i\le j]
\] 
denotes the ring of polynomials in the entries $\mu_k$ and
$\sigma_{ij}$ of the mean vector $\mu$ and the covariance matrix
$\Sigma$.

\begin{definition}
  A Gaussian model is an {\em implicit algebraic model} if its mean
  parameter space $M$ is equal to the intersection of an algebraic variety
  $V$ and the Cartesian product $\mathbb{R}^p\times \mathbb{R}^{p\times
    p}_{\rm pd}$. In other words, there exist polynomials
  $f_1,\dots,f_t\in\mathbb{R}[\mu_k,\sigma_{ij}\mid i,j,k\in[p],\; i\le j]$
  such that
  \[
  M = \left\{ (\mu,\Sigma)\in \mathbb{R}^p\times \mathbb{R}^{p\times
      p}_{\rm pd} \mid f_1(\mu,\Sigma)=\dots=f_t(\mu,\Sigma)=0
  \right\}.
  \]  
\end{definition}
The next fact is a consequence of the Tarski-Seidenberg Theorem
\cite[Thm.~2.3.4]{benedetti}.
\begin{proposition}
  \label{prop:semialg}
  The mean parameter space of an algebraic Gaussian model, parametric
  or implicit, is a semi-algebraic set.
\end{proposition}

An immediate consequence of Proposition \ref{prop:semialg} is that an
algebraic Gaussian model always has a well-defined dimension, namely,
the dimension $\dim(M)$ of the semi-algebraic mean parameter space
\cite[Def.~2.5.3]{benedetti}.  In the parametric case,
$\dim(M)=\dim(\mathbf{f}(\Theta))$ can be determined as the maximal
rank of the Jacobian of the rational parameterization map
$\mathbf{f}$.  In the implicit case, $\dim(M)$ can be computed based
on Gr\"obner basis techniques \cite[Thm.~3.7]{ascb}.
An implicit algebraic model need not be a parametric algebraic model, and
vice versa.  Nevertheless, a technique known as {\em implicitization\/},
cf. \cite[\S3.3]{coxlittleoshea} and \cite[\S3.2]{ascb}, permits to find a
unique smallest implicit model whose mean parameter space $M$ contains the
mean parameter space $\mathbf{f}(\Theta)$ of a given parametric model while
satisfying $\dim(M)=\dim(\mathbf{f}(\Theta))$.

\section{Identifiability}
\label{sec:ident}

When specifying a parametric statistical model, one of the first concerns
is whether the model is identifiable, that is, whether the parameters
uniquely specify probability distributions in the model.

\begin{definition}
  \label{def:ident}
  Consider a parametric Gaussian model with mean parameter space
  $M=\mathbf{f}(\Theta)$ given as the image of a map $\mathbf{f}:\Theta \to
  \mathbb{R}^p\times \mathbb{R}^{p\times p}_{\rm pd}$.  The model is
  \begin{enumerate}
  \item[(i)] {\em globally identifiable} at $\theta_0\in\Theta$ if 
    $\mathbf{f}^{-1}(\mathbf{f}(\theta_0))=\{\theta_0\}$;
  \item[(ii)] {\em locally identifiable} at $\theta_0\in\Theta$ if there
    exists a ball $B_\varepsilon(\theta_0)$ with center $\theta_0$ and
    radius $\varepsilon>0$ such that
    $\mathbf{f}^{-1}(\mathbf{f}(\theta_0))\cap B_\varepsilon(\theta_0)
    =\{\theta_0\}$.
  \end{enumerate}
  If the model is globally identifiable at all points in $\Theta$, that is,
  if the map $\mathbf{f}:\Theta \to M$ is a bijection, then we say that the
  model is {\em identifiable}.
\end{definition}

For parametric algebraic models, global and local identifiability at a
given point $\theta_0\in\Theta$ can be investigated by studying
whether a system of polynomial equations deduced from the possibly
rational equation system $\mathbf{f}(\theta_0)=\mathbf{f}(\theta)$ has
a (locally) unique solution $\theta\in\Theta$.  We illustrate this in
the following example.

\begin{figure}[t]
  \begin{center}
    \small
    \newcommand{\myNode}[2]{\circlenode{#1}{\makebox[2.75ex]{#2}}}
    \begin{pspicture}(-2,0)(2,2)
    \psset{linewidth=0.6pt,arrowscale=1.5 2,arrowinset=0.1}  
    \rput(0, 0.75){\myNode{h}{$H$}} 
    \psset{fillcolor=lightgray, fillstyle=solid}          
    \rput(-2, 1.5){\myNode{1}{$X_1$}}
    \rput(2, 1.5){\myNode{2}{$X_2$}} 
    \rput(-2, 0){\myNode{3}{$X_3$}} 
    \rput(2, 0){\myNode{4}{$X_4$}} 
    \ncline{->}{1}{h}
    \ncline{->}{2}{h}
    \ncline{->}{h}{3}
    \ncline{->}{h}{4}
    \end{pspicture}
  \end{center}
\caption{\label{fig:identgraph} Acyclic digraph for a hidden variable model.}
\end{figure}
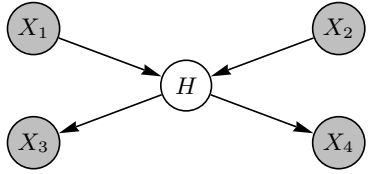

Consider the directed graphical Gaussian model with one hidden variable
depicted in Figure \ref{fig:identgraph}, which shows the relationship
between $p=4$ observed variables (shaded nodes) and one hidden variable
$H$.  For simplicity we consider the model comprising only centered
distributions.  This model is a parametric algebraic model with mean
parameter space $M=\mathbf{f}(\Theta)$ given by a polynomial map.  The
points in the parameterization domain $\Theta=\mathbb{R}^4\times
(0,\infty)^4$ are vectors
$\theta=(\beta_1,\beta_2,\beta_3,\beta_4,\omega_1,\omega_2,\omega_3,\omega_4)$.
Here, the four regression coefficients $\beta_i\in\mathbb{R}$ appear in
conditional means, namely $\mathbb{E}[H\mid X_1,X_2]=\beta_1X_1+\beta_2X_2$
and $\mathbb{E}[X_i\mid H]=\beta_iH$ if $i=3,4$.  The variances
$\omega_i>0$ are either marginal or conditional variances, $\omega_i
=\text{Var}[X_i]$ for $i=1,2$, and $\omega_i =\text{Var}[X_i\mid H]$ for
$i=3,4$.  The parameterization map is
\[
\begin{split}
  \mathbf{f}: \Theta &\to
  \mathbb{R}^4\times\mathbb{R}^{4\times 4}_{\rm pd}\\
  \theta
  &\mapsto 
  [ 0,\mathbf{f}_\Sigma(\theta)],
\end{split}
\]
where $\mathbf{f}_\Sigma(\theta)$ is the symmetric covariance matrix 
\begin{equation}
\label{eq:identmap}
  \begin{pmatrix}
    \omega_1& 0& \beta_3\beta_1\omega_1& \beta_4\beta_1\omega_1\\ 
    & \omega_2 & \beta_3\beta_2\omega_2 & \beta_4\beta_2\omega_2\\
     & 
     & 
     \omega_3+\beta_3^2(\beta_1^2\omega_1 + \beta_2^2\omega_2 + 1)  &
    \beta_4\beta_3(\beta_1^2\omega_1 + \beta_2^2\omega_2 + 1)\\ 
     &  &
    &
    \omega_4+\beta_4^2(\beta_1^2\omega_1 + \beta_2^2\omega_2 + 1)
  \end{pmatrix}.
\end{equation}
Note that we set the conditional variance $\text{Var}[H\mid X_1,X_2]=1$
because the image of the parameterization map remains unchanged if a free
parameter for this variance is introduced. For details on such
parameterizations see e.g.~\cite[\S8]{rich:2002}.

Is this parametric hidden variable model globally (or locally)
identifiable at $\theta_0\in\Theta$?  We can answer this question by
studying the system of equations $\mathbf{f}_\Sigma(\theta_0) =
\mathbf{f}_\Sigma(\theta)$ in which the components $\beta_{i0}$ and
$\omega_{i0}$ of $\theta_0$ are fixed numbers and the components $\beta_i$
and $\omega_i$ of $\theta$ are indeterminants.  From (\ref{eq:identmap}),
it is apparent that if $\mathbf{f}(\theta) = \mathbf{f}(\theta_0)$ then
$\omega_1=\omega_{10}$ and $\omega_2=\omega_{20}$.  Additional consequences
can be worked out by hand, but we can also let the computer do this for us.

\begin{table}[!t]
\begin{small}
\begin{verbatim}
LIB "linalg.lib";  option(redSB);
ring R = (0,b10,b20,b30,b40,w10,w20,w30,w40),
         (b1,b2,b3,b4,w1,w2,w3,w4),dp;
// b1,...,w4 are indeterminants; b10,...,w40 are symbolic parameters
matrix B[5][5] = 1,0,0,0,0,
                 0,1,0,0,0,
                 0,0,1,0,-b3,
                 0,0,0,1,-b4,
                 -b1,-b2,0,0,1;
matrix W[5][5] = w1,0,0,0,0,
                 0,w2,0,0,0,
                 0,0,w3,0,0,
                 0,0,0,w4,0,
                 0,0,0,0,1;
matrix B0[5][5] = 1,0,0,0,0,
                 0,1,0,0,0,
                 0,0,1,0,-b30,
                 0,0,0,1,-b40,
                 -b10,-b20,0,0,1;
matrix W0[5][5] = w10,0,0,0,0,
                 0,w20,0,0,0,
                 0,0,w30,0,0,
                 0,0,0,w40,0,
                 0,0,0,0,1;
matrix f[4][4] = submat(inverse(B)*W*inverse(transpose(B)),1..4,1..4);
matrix f0[4][4] = submat(inverse(B0)*W0*inverse(transpose(B0)),1..4,1..4);
ideal I=0;  int i,j;
for(i=1; i<=4; i++){ for(j=i; j<=4; j++){
    I = I + ideal(f0[i,j]-f[i,j]); // identifiability equations
} }
ideal J = std(I);                  // Groebner basis for ideal I
dim(J); mult(J);  
\end{verbatim}
\end{small}
\caption{\label{tab:sing} Code from session in computer algebra system {\tt
  Singular}.}
\end{table}

Running the code in Table \ref{tab:sing} in the computer algebra system
{\tt Singular} \cite{greuel} informs us that the model is not globally
identifiable at generic $\theta_0\in\Theta$ because the solution set
$\mathbf{f}^{-1}(\mathbf{f}(\theta_0))$ generally contains
$\mathtt{mult(J)}=2$ isolated points ($\mathtt{dim(J)}=0$).  The computed
{\em Gr\"obner basis\/} (see \cite{coxlittleoshea,ascb} for the relevant
background)
\begin{verbatim}
> J;
J[1]=w4+(-w40)
J[2]=w3+(-w30)
J[3]=w2+(-w20)
J[4]=w1+(-w10)
J[5]=(b40)*b3+(-b30)*b4
J[6]=(-b40)*b2+(b20)*b4
J[7]=(-b40)*b1+(b10)*b4
J[8]=b4^2+(-b40^2)
\end{verbatim}
suggests that for $\beta_{40}\not=0$, it holds that
\begin{equation}
  \label{eq:ident02}
  \mathbf{f}^{-1}(\mathbf{f}(\theta_0)) = \{
  (\beta_{10},\dots,\beta_{40},\omega_{10},\dots,\omega_{40})^t,
  (-\beta_{10},\dots,-\beta_{40},\omega_{10},\dots,\omega_{40})^t
  \}.
\end{equation}
However, in the Gr\"obner basis computation simplifications are made that
are only valid if {\tt b10,...,w40} are generic.  In other words, during
the computation a (finite) number of polynomial expressions in {\tt
  b10,...,w40} are assumed to be non-zero.  So while we can conclude that
(\ref{eq:ident02}) holds for almost every $\theta_0\in\Theta$, it may and
does fail at certain points in the parameter domain $\Theta$.  For example,
(\ref{eq:ident02}) does not hold if $\beta_{40}=0$, in which case it is
possible that $\mathtt{dim(J)}\in\{ 1,2\}$ indicating failure of local
identifiability.  In conclusion, the computation shows that the model is
locally identifiable at almost every $\theta_0\in\Theta$.
In more complicated models, computations treating $\theta_0$ as
symbolic quantity may become prohibitive.  However, solving the system
$\mathbf{f}(\theta_0)=\mathbf{f}(\theta)$ for a particular numeric
vector $\theta_0$ may still be feasible and informative.


\section{Singularities of Gaussian models}
\label{sec:sing}

The saturated Gaussian model is a regular exponential family
(Proposition~\ref{prop:expfam}).  Therefore, under ``regularity''
conditions, results about asymptotic distributions of MLE and likelihood
ratio test statistic can be transfered to submodels.  Suppose the submodel
is an algebraic model with mean parameter space $M$ containing the true
distribution $\mathcal{N}(\mu_0,\Sigma_0)$.  If $M$ is a smooth manifold,
in which case the submodel is a {\em curved exponential family}, then
regardless of where the true parameter $(\mu_0,\Sigma_0)$ is located, the
MLE is asymptotically normal as the sample size tends to infinity.
Moreover, the likelihood ratio test statistic for testing the submodel
against the saturated model has an asymptotic $\chi^2$-distribution with
degrees of freedom equal to the codimension of $M$, that is, the difference
between the dimension of the saturated mean parameter space and $\dim(M)$.

These standard results need no longer be true if one leaves the
realm of curved exponential families.  For example, if inequality
constraints are imposed on the mean parameter space of a curved exponential
family, boundary effects may be created.  More subtly, the regularity
conditions may be violated at points that are ``singularities'' in the mean
parameter space of an algebraic Gaussian model.  For a rigorous definition
of singularities of algebraic varieties, see e.g.~\cite[\S3.2]{benedetti}.
In the examples in this section the singularities are obvious and
intuitive.  However, this will not necessarily be the case in larger
models, in which case computer algebra software is very helpful for
locating singular points.  In particular, the software {\tt Singular}
offers the command {\tt slocus} for computation of singular loci.

\subsection{Simple bivariate examples under independence}
\label{subsec:euclidgeo}

Issues with singularities can be illustrated nicely with bivariate normal
distributions.  For a closed set $C\subseteq \mathbb{R}^2$, let $M_C
=C\times \{I_2\}$ be the mean (and natural) parameter space of the model of
all bivariate normal distributions with mean vector $\mu\in C$ and
covariance matrix $\Sigma$ equal to the identity matrix
$I_2\in\mathbb{R}^{2\times 2}$.  In this case the MLE $\hat\mu$ for the
model with mean parameter space $M_C$ is the point in $C$ that is closest
to the sample mean vector $\bar X$ in Euclidean distance.  The likelihood
ratio test statistic $\lambda_C$ for testing $\mu\in M_C$ versus $\mu\in
M_{\mathbb{R}^2}$ is equal to the product of the sample size $n$ and the
squared Euclidean distance of $\bar X$ and $C$.


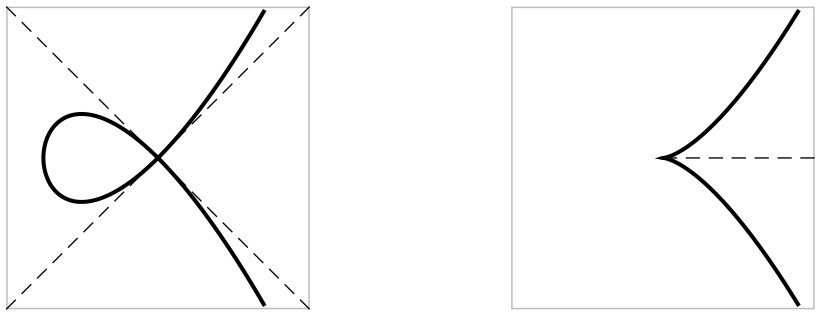
\begin{figure}[t]
  \centering
  \begin{tabular}{ll}
\hspace{1.5cm}
\begin{minipage}{6cm}
\begin{pspicture}(-1.325,-1.325)(1.325,1.95)
\psset{xunit=1.5cm,yunit=1.5cm}
\psframe[linecolor=lightgray,linewidth=.2mm](-1.325,-1.325)(1.325,1.325)
\parametricplot[plotstyle=curve,plotpoints=200,%
linewidth=1.5pt]{-1.39}{1.39}%
{t t mul 1 sub t t mul 1 sub t mul}
\psline[linecolor=black,linestyle=dashed,linewidth=.5pt](-1.325,-1.325)(1.325,1.325)
\psline[linecolor=black,linestyle=dashed,linewidth=.5pt](1.325,-1.325)(-1.325,1.325)
\end{pspicture}
\vspace{.21cm}
\end{minipage}
&
\hspace{-.1cm}
\begin{minipage}{6cm}
\begin{pspicture}(-1.5,-1.325)(1.325,1.95)
\psset{xunit=1.5cm,yunit=1.5cm}
\psframe[linecolor=lightgray,linewidth=.2mm](-1.325,-1.325)(1.325,1.325)
\parametricplot[plotstyle=curve,plotpoints=200,%
linewidth=1.5pt]{-1.09}{1.09}%
{t t mul t t t mul mul}
\psline[linestyle=dashed,linewidth=.5pt](0,0)(1.325,0)
\end{pspicture}
\vspace{.21cm}
\end{minipage}
  \end{tabular}
  \vspace{0.45cm}
  \caption{The parameter spaces of two simple algebraic Gaussian models.}
  \label{fig:3simple}  
\end{figure}

\begin{example}[Folium of Descartes]
  \label{ex:folium}
  \rm Let $C=\{\mu\in\mathbb{R}^2\mid \mu_2^2=\mu_1^3+\mu_1^2\}$, which is
  a curve that can be parameterized as
  $\mathbf{f}(\theta)=[\theta^2-1,\theta(\theta^2-1)]$.  The curve is shown
  in the left plot in Figure \ref{fig:3simple}. The algebraic model with
  mean parameter space $M_C$ is {\em not\/} a curved exponential family due
  to the singularity at the point of self-intersection, which is $\mu=0$.
  The dashed lines $\mu_2=\pm\mu_1$ in the plot are orthogonal to each
  other and indicate the tangent cone at $\mu=0$.  If the true parameter
  point is $\mu=0$, then the asymptotic distribution of the likelihood
  ratio test statistic $\lambda_C$ is given by the squared Euclidean
  distance between a draw from $\mathcal{N}(0,I_2)$ and the two orthogonal
  lines.  This asymptotic distribution is the distribution of the minimum
  of two independent $\chi^2_1$-random variables.  \qed
\end{example}

\begin{example}[Neil's parabola]
  \label{ex:neil}
  \rm Let $C=\{\mu\in\mathbb{R}^2\mid \mu_2^2=\mu_1^3\}$ be the curve with
  parameterization $\mathbf{f}(\theta)=(\theta^2,\theta^3)$, which is shown
  in the right-most picture of Figure \ref{fig:3simple}.  The algebraic
  model with mean parameter space $M_C$ is again {\em not\/} a curved
  exponential family due to the singularity at the cusp point $\mu=0$.  For
  true parameter point $\mu=0$, the likelihood ratio test statistic
  $\lambda_C$ has an asymptotic distribution that is the mixture of a
  $\chi^2_1$- and a $\chi^2_2$-distribution.  This mixture distribution is
  the distribution of the squared Euclidean distance between a draw from
  $\mathcal{N}(0,I_2)$ and the (dashed) half-ray $\{\mu\mid \mu_1\ge
  0,\mu_2=0\}$.  \qed
\end{example}

Examples \ref{ex:folium} and \ref{ex:neil} demonstrate non-standard
asymptotics at model singularities.  At regular points in the
respective mean parameter spaces the usual asymptotics apply.
However, if the true parameter forms a regular point that is close to
the singular locus then a very large sample size may be required in
order for the usual asymptotics to provide good approximations to the
distributions of estimators and test statistics.


\subsection{A conditional independence model with singularities}
\label{subsec:an-indep-model}

Many conditional independence models, in particular graphical models, form
curved exponential families.  However, singularities may arise from
combining arbitrary independence constraints.  Consider, for example, the
model of trivariate normal distributions under which a random vector
satisfies $X_1\ind X_2$ and simultaneously $X_1\ind X_2\mid X_3$.  By
Proposition~\ref{prop:ci}, the model is an implicit algebraic model with
mean parameter space
\[
M = \left\{ (\mu,\Sigma)\in\mathbb{R}^3\times\mathbb{R}^{3\times
    3}_{\rm pd} \mid \sigma_{12}=0,\; \det(\Sigma_{\{1,3\}\times
    \{2,3\}})= \sigma_{12}\sigma_{33}-\sigma_{13}\sigma_{23}= 0
\right\}
\]
The set $M$ is defined by the joint vanishing of the two polynomials
$f_1=\sigma_{12}$ and $f_2=\sigma_{13}\sigma_{23}$.  We see that 
\[
M = M_{13}\cup M_{23} := 
\{(\mu,\Sigma)\in M\mid \sigma_{12}=\sigma_{13}=0\}
\cup
\{(\mu,\Sigma)\in M\mid \sigma_{12}=\sigma_{23}=0\}.
\] 
This reflects the well-known fact that
\[
\left[\;X_1\ind X_2\, \wedge\, X_1\ind X_2\mid X_3 \;\right] \iff
\left[\;X_1\ind (X_2,X_3)\, \vee\, X_2\ind (X_1,X_3) \;\right],
\]
which also holds for distributions other than the multivariate normal;
compare \cite[Thm.~8.3]{dawid:1980}.  The singular locus of $M$ is the
intersection
\[
M_{\rm sing} = M_{13}\cap M_{23} = \{(\mu,\Sigma)\in M\mid
\sigma_{12}=\sigma_{13}=\sigma_{23}=0\},
\]
which corresponds to diagonal $\Sigma$, i.e., complete
independence $X_1\ind X_2\ind X_3$.

The likelihood ratio test statistic for testing the model with mean
parameter space $M$ against the saturated model can be expressed as
\begin{equation}
  \label{eq:cilr}
\lambda = n\cdot \min\left\{ \log\left(
\frac{s_{11}\det(S_{\{2,3\}\times\{2,3\}})}{\det(S)}\right),
\log\left(
\frac{s_{22}\det(S_{\{1,3\}\times\{1,3\}})}{\det(S)}
\right)\right\}.
\end{equation}
If $(\mu,\Sigma)\in M\setminus M_{\rm sing}$, then $\lambda$ converges
to a $\chi^2_2$-distribution for $n\to\infty$; over the singular locus
the limiting distribution is non-standard.

\begin{proposition}
  \label{prop:cisingasy}
  Let $(\mu,\Sigma)\in M_{\rm sing}$, and let $W_{12}$, $W_{13}$,
  $W_{23}$ be three independent $\chi^2_1$-random variables.  As
  $n\to\infty$, the likelihood ratio test statistic $\lambda$
  converges to the minimum of two dependent $\chi^2_2$-distributed
  random variables, namely,
  \[
  \lambda \longrightarrow_d
  \min(W_{12}+W_{13},W_{12}+W_{23})=
  W_{12}+\min(W_{13},W_{23}).
  \]
\end{proposition}
\begin{proof}
   For $i\in\{1,2,3\}$ and $A\subset \{1,2,3\}$, let 
  \[
  s_{ii.A} = s_{ii}-S_{\{i\}\times A}S_{A\times A}^{-1}S_{A\times \{i\}}.
  \]
  The likelihood ratio test statistic can be rewritten as
  \begin{equation}
    \label{eq:cilr2}
    \lambda = n \log\left(
        \frac{s_{11}s_{22}}{s_{11}s_{22}-s_{12}^2}
        \right)
        +
        n\cdot \min\left\{ \log\left(
        \frac{s_{33.2}}{s_{33.12}}\right),
      \log\left(
        \frac{s_{33.1}}{s_{33.12}}
      \right)\right\}.
  \end{equation}
  Recall that 
  \begin{multline}
    \label{eq:ciiss}
    \sqrt{n}\big[(s_{11},s_{12},s_{13},s_{22},s_{23},s_{33})^t-
    (\sigma_{11},\sigma_{12},\sigma_{13},
  \sigma_{22},\sigma_{23},\sigma_{33})^t
  \big]
  \longrightarrow_d \\
  \mathcal{N}\left(0,
    (\sigma_{ik}\sigma_{jm}+\sigma_{im}\sigma_{jk})_{ij,km}\right).
  \end{multline}
  Since $(\mu,\Sigma)\in M_{\rm sing}$ implies that $\Sigma$ is diagonal,
  the covariance matrix of the normal distribution in (\ref{eq:ciiss}),
  known as the Isserlis matrix of $\Sigma$, is diagonal.  Using an
  expansion up to the Hessian in the delta-method
  \cite[\S3.3]{vandervaart}, we can show that the three logarithmic terms
  in (\ref{eq:cilr2}) converge to three independent $\chi^2_1$-random
  variables.
\end{proof}

\section{Conclusion}
\label{sec:conc}

The goal of this paper was to demonstrate the usefulness of algebraic
geometry for studying properties of statistical (Gaussian) models.  In
order to keep intuition alive, the examples in this paper were chosen to be
rather simple, but algebraic geometry can also provide useful insights in
larger, less tractable models.  

Two particular problems were visited in this paper.  First, parameter
identifiability often gives rise to polynomial equation systems, the
structure of which becomes more transparent when the equations are
presented in Gr\"obner basis form (Section \ref{sec:ident}).  Second, model
singularities can result into non-standard asymptotics (Section
\ref{sec:sing}).  Locating singularities and working out the associated
asymptotics are the first steps towards solving the challenging problem of
divising sensible statistical procedures for models with singularities.
Finally, we remark that methods combining Gr\"obner basis techniques with
numerical solving can also be used to compute all solutions to interesting
likelihood equations, compare e.g.~\cite{drton:06}.

\authorsPS
\end{document}